\documentclass{amsart}
\usepackage{amscd}
\usepackage{verbatim}
\usepackage{amssymb}
\usepackage{epic,eepic}

\usepackage{epsfig}

\newcommand\Mand{\ \text{and}\ }

\newcommand\Mor{\ \text{or}\ }

\newcommand\Id{\operatorname{Id}}
\newcommand\Real{\mathbb{R}}

\newcommand\im{\operatorname{Im}}

\newcommand\pa{\partial}
\newcommand\Rn{\Real^n}

\newcommand\Cinf{{\mathcal C}^{\infty}}
\newcommand\dist{{\mathcal C}^{-\infty}}

\newcommand\supp{\operatorname{supp}}

\newcommand\bop{{\mathcal B}}

\newcommand\sci{{}^{\text{sc}}}
\newcommand\sct{\sci T^*}

\newcommand\scct{\sci\bar{T}^*}

\newcommand\Hsc{H_{\text{sc}}}

\newcommand\Psop{\operatorname{\Psi}}

\newcommand\psit{\tilde\psi}

\newcommand\bX{\partial X}

\newcommand\bl{{\text b}}
\newcommand\scl{{\text{sc}}}
\newcommand\sccl{{\text{scc}}}

\newcommand\Psisch{\Psop_{\scl,h}}
\newcommand\Psiscch{\Psop_{\sccl,h}}

\newcommand\Vb{{\mathcal V}_{\bl}}

\newcommand\ep{\epsilon}

\newtheorem{lemma}{Lemma}[section]
\newtheorem{prop}[lemma]{Proposition}

\newtheorem*{thm*}{Theorem}
\newtheorem*{prop*}{Proposition}
\newtheorem*{conj*}{Conjecture}
\numberwithin{equation}{section}
\theoremstyle{remark}
\newtheorem{rem}[lemma]{Remark}
\theoremstyle{definition}

\newtheorem*{Def*}{Definition}

\title{Semiclassical estimates in asymptotically Euclidean scattering}
\author{Andr\'as Vasy}
\author{Maciej Zworski}
\address{Mathematics Department, University of California, Berkeley}
\email{andras@math.berkeley.edu}
\email{zworski@math.berkeley.edu}
\date{October 19, 1999}

\begin{document}
\maketitle

\section{Introduction}
The purpose of this note is to obtain semiclassical resolvent
estimates for long range perturbations of the Laplacian on 
asymptotically Euclidean manifolds.

For an estimate which is uniform in the Planck constant $ h$
we need to assume that the energy level is  {\em non-trapping}.
In the high energy limit (that is, when we consider $ \Delta - \lambda^2 $, 
as
$ \lambda \rightarrow \infty $, which is equivalent to $ h^2 \Delta 
- 1 $, $ h \rightarrow 0 $), this corresponds to the global assumption
that the geodesic flow is non-trapping. We note here that
a sufficiently small neighbourhood of infinity is always non-trapping.

The resolvent estimates in the classical ($h=1$) and semi-classical
cases have a long tradition going back to the {\em limiting absorption
principle} -- see \cite{Ag} and references given there. Various
variants of the 
theorem we present were proved in Euclidean potential scattering 
by Jensen-Mourre-Perry \cite{JeMoPe}, Robert-Tamura \cite{RoTa}, 
G\'erard-Martinez \cite{GeMa}, G\'erard \cite{Ge} and Wang \cite{Wa}.
The proofs were based on {\em Mourre theory} whose underlying feature
is the positive commutator method accompanied by functional analytic
techniques for obtaining a resolvent estimate.
While the work of G\'erard-Martinez \cite{GeMa} explains the
role of geometry in the positive commutator estimate itself, it
refers to Mourre's work for the
functional analytic argument. We adopt a completely geometric approach
based on direct microlocal ideas.

The classical version of the 
estimate on asymptotically Euclidean manifolds
($h=1$ in which case there is no need for the non-trapping
assumption) is essentially in Melrose's
original paper on the subject \cite{mel:spectral} in which he introduced
a fully microlocal point of view to scattering. 
However, the proof presented here is somewhat
different in spirit: a global positive commutator argument is used to
derive an estimate on the resolvent directly.

Referring to \eqref{eq:scatmet} below for the definition of a scattering
metric, to \eqref{eq:pd},\eqref{eq:vd} for the definition of a 
long range semi-classical perturbation, and \eqref{eq:non-trapping}
for the definition of a non-trapping energy, we state our main

\medskip

\noindent
{\bf Theorem.} {\em Let $ X $ be a manifold with boundary and let
$ \Delta $ be the Laplacian of a scattering metric on $ X$. If 
$ P = h^2 \Delta + V  $ is a semi-classical
long range perturbation
of $ h^2 \Delta $,
and $ R ( \lambda ) = ( P - \lambda )^{-1} $ its
resolvent, then for all $m\in\Real$,
\begin{equation}
\|R(\lambda+it)f\|_{\Hsc^{m,-1/2-\ep}(X)}\leq C_0h^{-1}\|f\|
_{\Hsc^{m-2,1/2+\ep}(X)},\quad\ep>0,
\end{equation}
with $C_0$ independent of $t\neq 0$ real and $\lambda\in I$,
$I\subset(0,+\infty)$ a
compact interval in the set of non-trapping energies for $P$.}

\medskip

Here $ \Hsc^{m,k}(X) $ denote Sobolev spaces adapted to the scattering
calculus, that is to asymptotically Euclidean structures. The index $ m$
indicates smoothness and $ k$ the rate of decay at infinity: the larger
the better in both cases.

To indicate the main idea of the proof let
$p-\lambda$ be the principal symbol of
$P-\lambda$. Here the principal symbol is meant in both the semi-classical
sense and the scattering sense -- see Sect.2. 
Near the characteristic variety of $P-\lambda$, we construct a function
$q\geq 0$ such that $q$ is decreasing along the Hamilton vector field
$H_p$. This gives the required estimate for the resolvent when we apply 
a variant of the well known commutator method -- see \cite{Ho:e} for the
now standard application to the propagation of singularities for operators of
principal type. 

We stress that to prove, say, the outgoing resolvent estimate, one needs to keep
the signs of both $q$ and $H_p q$ fixed throughout phase space, and in
case of the outgoing estimate, these signs must be opposite. Indeed,
it is the fixed sign of $q$ that makes it possible to eliminate the machinery
of Mourre's method. The positivity
of $q$ shows that in the outgoing region, where bicharacteristics
tend as $t\to+\infty$, $q$ must be of the form $x^r a$, $a\in\Cinf(\sct X)$,
$r>0$, (here $x$ is a boundary defining function), and in the incoming
region, where bicharacteristics tend as $t\to-\infty$, it must be of
the form $x^{-s} b$, $b\in\Cinf(\sct X)$, $s>0$. The difference between
these two weights, which can be made arbitrarily small, but is never $0$,
plus the improvement by $1$ in the order when calculating a commutator,
explains how the weighting of the Sobolev spaces works.

For applications of the non-trapping estimates to more general operators
we refer to a recent paper by Bruneau-Petkov \cite{BrPe}. It is clear 
that the ``black box'' set-up discussed there can be easily adapted to 
the manifold situation.

\medskip

\noindent
{\sc Acknowledgments.}
The authors are grateful to the National Science Foundation for partial support.

\section{Preliminaries}

Let $X$ be a ${ \mathcal C}^\infty$
 manifold with boundary, $ \partial X  $ and let 
$x$ be a  boundary defining function. 
Thus, in a small collar neighborhood $[0,\ep_0)\times\bX$
of the boundary $\bX$, we have `semi-global
coordinates' $(x,y,\xi,\eta)$ on $ T^*_{[0,\ep_0)\times\bX} X $.

Microlocal techniques adapted to asymptotically Euclidean structure
near $ \bX $ (see \eqref{eq:scatmet})
were introduced by Melrose \cite{mel:spectral}. We
start by recalling the scattering cotangent bundle $ \sct X $ which is
the natural {\em phase space}.
It is defined as the dual 
of the scattering tangent bundle $ \sci T X$,  
which in turn is defined so that
the space of vector fields $ x {\mathcal V}_b ( X) $, 
where $ {\mathcal V}_b ( X) $ 
are vector fields tangent to $\bX $, is given by sections
\[ x {\mathcal V}_b ( X) = {\mathcal C}^\infty ( X ; \sci  T X) \,,\]
see \cite{mel:spectral} for a thorough discussion. 
Since $ \sci T X \hookrightarrow   T X $ we have a natural
map $  T^* X \rightarrow \sct X $. In `semi-global
coordinates' $(x,y,\xi,\eta)$ on $ \sct_{[0,\ep_0)\times\bX} X $
it is given by 
\[ ( x, y , \tau, \mu ) = ( x, y , x^2 \xi , x \eta ) \,, \]
and this identification is worth keeping in mind since the symplectic
and contact structures are inherited from $ T^* X $, that is, from 
the $ ( x, \xi ) $ coordinates. In particular, when we speak of
the Hamilton vector fields on $ \sct X $, we mean the natural 
extention of the usual Hamilton vector field on $ \sct X^o \simeq
T^* X^o $, to $ \sct X$ -- see \cite{mel:spectral} and \cite{MelZw}.
We also note that the variable $ \mu $ is naturally identified with 
$\mu\in T^*_y\bX$. 

The fiber radial compactification of $\sct X$ is
denoted by $\scct X$; $\scct X$ is thus a ball bundle over $X$.
Classical symbols, $a\in S^{m,l}_{cl}(X)$,
are functions $a\in x^{l}\rho_\infty^{-m} \Cinf(\scct X)$.
By $a\in S^{m,l}(X)$ it is meant that $a\in\Cinf(T^*X)$,
$x^{-l}\rho_\infty^m a\in L^\infty(T^*X)$,
and the same estimate holds after the application (to $a$)
of any $b$-differential
operator on $\scct X$, that is, an operator in the algebra 
generated by $ {\mathcal V}_b ( \scct X ) $, vector fields tangent to 
$ \partial \scct X $.  

\begin{figure}[ht]
\begin{center}
\mbox{\epsfig{file=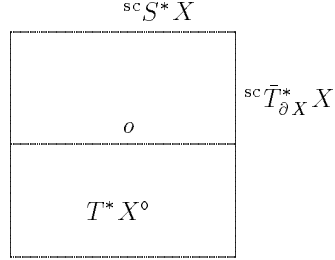}}
\end{center}
\caption{The fiber compactification $\scct X$ of $\sct X$ is a manifold
with corners. Its boundary hyperfaces are $\scct_{\bX}X$, which is
a ball bundle over $\bX$, and the cosphere bundle ${}^{\text{sc}}S^* X$,
which is
a sphere bundle over $X$. The zero section is denoted by $o$.}
\label{fig:scct}
\end{figure}

The semiclassical calculus for the Weyl metric on $ T^* X$.
\[ \frac{dz^2}{1 + |z|^2 } + \frac{d\zeta^2}{ 1 + |\zeta|^2 } \]
is well known and, for instance, it is discussed in great generality 
in \cite{HelSj}. The natural generalization to manifolds with 
asymptotically Euclidean structure near infinity is given in the 
Appendix to \cite{WZ1}. We will review and slightly extend it below.

The semi-classical symbols are defined as follows:
 $a\in S^{m,l,k}(X)$ means that $a\in\Cinf((0,1)\times T^*X)$,
$h^k x^{-l}\rho_\infty^m a\in L^\infty((0,1)\times T^*X)$, and
the same estimate holds after the application of any $b$-differential
operator $ \scct X$. Thus, $a(h,.)\in S^{m,l}(X)$ for $h\in(0,1)$, and the
symbol estimates are uniform in $h$.
The corresponding {\large class} of classical
symbols, $a\in S^{m,l,k}_{cl}(X)$ are functions with 
$h^k x^{-l}\rho_\infty^m a\in\Cinf([0,1)\times \scct X)$.

\begin{figure}[ht]
\begin{center}
\mbox{\epsfig{file=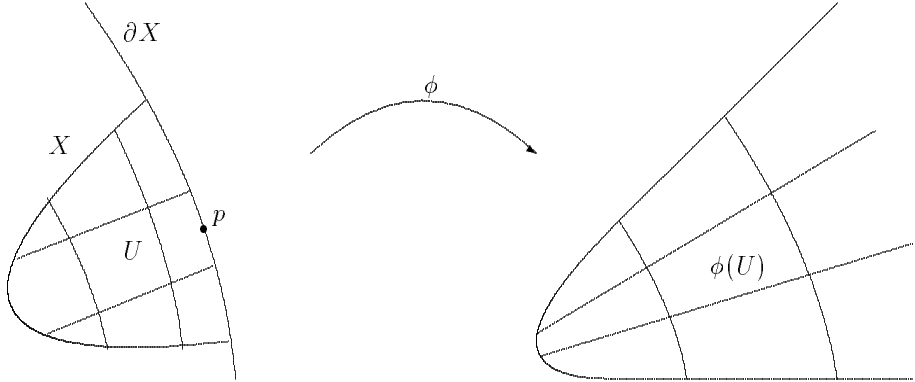}}
\end{center}
\caption{Local Euclidean coordinates near $p\in\bX$ identify a neighborhood
$U$ of $p$ in $X$ with a conic neighborhood $\phi(U)$ of infinity in $\Rn$.}
\label{fig:conic}
\end{figure}

For $ a \in S^{m,l,k} ( X) $ we define a semiclassical operator
$ \mathrm{Op} ( a) \in \Psi_{\mathrm{sc}}^{m,l,k} (X)$ 
as in Appendix to \cite{WZ1}:
we first use local Euclidean coordinates in a cone near infinity, 
identified
with a neighbourhood of a boundary point (see Figure~\ref{fig:conic})
to define
$$
  A u(z) = \left( \frac 1{2\pi h} \right)^n \int e^{i (z-w) 
\cdot \xi/h} a(h, z,\xi)
  u(w) dw d\xi.
  $$
with $a \in S^{m,l,k} ({\mathbb R}^n)$. Invariance
under local changes of coordinates then gives $ \mathrm{Op} (a) $ and
leads to the definition of the class $ \Psi_{\mathrm{sc}}^{m,l,k} (X)$.

We then have the symbol map $ \sigma_{\mathrm{sc},h} : 
 \Psi_{h, \mathrm{sc}}^{m,l,k} (X) \rightarrow S^{m,l,k}(X)  $ with the
usual properties, and in particular with the short exact sequence
\[ 
0 \to  \Psi_{h, \mathrm{sc}}^{m-1,l+1,k-1} (X) \to
\Psi_{\mathrm{sc}}^{m,l,k} (X) 
\stackrel{\sigma_{h, \mathrm{sc}}^{m,l,k}}\longrightarrow S^{m,l,k} ( X) / S^{m-1,l+1,k-1}
(X) \longrightarrow 0 \,.\]

Another important property of $\Psiscch(X)$ is that it is
commutative to top order, and the principal symbol of a commutator
is given by the Poisson bracket of the principal symbols of the commutants.
That is, if $A\in\Psiscch^{m,l,k}(X)$, $B\in\Psiscch^{m',l',k'}(X)$, then
$[A,B]\in\Psiscch^{m+m'-1,l+l'+1,k+k'-1}(X)$ and
\begin{equation}
\sigma_{h, \mathrm{sc}}^{m+m'-1,l+l'+1,k+k'-1}([A,B])=\frac hi H_a b,
\end{equation}
where $a$, $b$ are the principal symbols of $A$ and $B$, and
$H_a$ denotes the Hamilton vector field of $a$.

We will also make use of the sharp G\aa rding estimate:

\begin{lemma}
Suppose $A\in\Psiscch^{0,0,0}(X)$ is self-adjoint,
and its (joint semiclassical) principal
symbol is $a\geq 0$. Then there exists $C>0$ such that
\begin{equation}
\langle u,Au\rangle \geq -Ch\|u\|_{\Hsc^{-1/2,1/2}(X)}.
\end{equation}
In particular, 
if $A\in\Psiscch^{2m,-2l,0}(X)$ has principal symbol $a\geq 0$,
then $A\geq hR$ for some
$R\in\Psiscch^{2m-1,-2l+1,0}(X)$.
\end{lemma}
\begin{proof}
The inequality is well known in the case of $ {\mathbb R}^n $ -- 
see Sect.18.4 of \cite{ho:v3} (easily adapted to the semi-classical
setting), \cite{disj} and \cite{HelSj}. The localization argument 
presented in the Appendix of \cite{WZ1} then gives the lemma.
\end{proof}

Now, let $ g $ be a scattering metric on $ X $, that is, a metric
which near $ \bX $ takes the form
\begin{equation}
\label{eq:scatmet}
 \frac{dx^2}{x^4} + \frac{h'}{x^2} \,, \ \ h'|_{\bX} = h \ \text{is 
a metric on } \bX \,. \end{equation}
This defines an asymptotically Euclidean structure near $ \bX $:
a neighbourhood of $ \bX $ is isometric to a perturbation of the large
end of the cone $ {\mathbb R}_+ \times \bX $ with the metric $ 
dr^2 + r^2 h $. 

We will consider the following self-adjoint, classically elliptic 
operators in $ {\mathrm{Diff}_{h, \mathrm{sc}}}^{2,0,0}
( X )\subset \Psi_{h, \mathrm{sc}}^{2,0,0} $:
\begin{equation}
\label{eq:pd}
P = h^2 \Delta_g + V 
\end{equation}
where in any compact set, $ V $ is a second order semiclassical 
operator ($ V = \sum_{|\alpha| \leq 2} v_\alpha ( z , h ) ( hD_z)^\alpha
$ in local coordinates) and near the boundary $ \bX $, in local 
coordinates $ y \in \bX $, 
\begin{equation}\begin{split}
\label{eq:vd}
V = x^\gamma V_0 \,,& \ \ V_0 = \sum_{ |\alpha | + k \leq 2}
v_{ k \alpha } ( x, y , h ) (h x^2 D_x)^k (hD_y)^\alpha \,,\\
&v_{ k \alpha }-v^0_{k\alpha} \in hS^{0,0,0} ( X)\,,\ v^0_{k\alpha}\in
S^{0,0}(X) \ \ \gamma > 0 \,.
\end{split}\end{equation}
The condition that the coefficients are symbols independent of the 
fiber variables means that $ |(x\partial_x)^l \partial^\beta_y 
v_{ k , \alpha } | \leq C_{ l \beta } $. In the Euclidean setting
it corresponds to assuming that the coefficients are symbols in 
the Euclidean base variables. Due to the vanishing of $v_{k\alpha}
-v^0_{k\alpha}$ in $S^{0,0,0}(X)$ when $h=0$, the semiclassical
principal symbol of $P$ is
\begin{equation}
p=g+x^\gamma \sum_{ |\alpha | + k \leq 2}v^0_{k\alpha}(x,y)\tau^k\mu^\alpha,
\end{equation}
where $g$ also denotes the (dual) metric function of the metric $g$. Thus,
$p$ can be represented by an $h$-independent function, which will
be convenient for the construction in the last section of this paper.
Note, however, that in \eqref{eq:vd},
$v_{ k \alpha }-v^0_{k\alpha} \in hS^{0,0,0} ( X)$ could be replaced by
$v_{ k \alpha }-v^0_{k\alpha} \in h^\rho S^{0,0,0} ( X)$, $\rho>0$,
or indeed by the assumption that $v_{k\alpha}$ is continuous on $[0,1)_h$
with values in $S^{0,0}(X)$, at the expense of minor changes in the next
section.

For obtaining the uniform resolvant estimates in $h$ for $R(\lambda\pm
i0)$, we make the
assumption that the Hamiltonian is non-trapping at energy $\lambda$
\begin{equation}\begin{split}\label{eq:non-trapping}
\text{for any}\ \xi\in T^* X^\circ
&\ \text{satisfying}\ p(\xi)=\lambda,\\
&\ \lim_{t\to\pm\infty}x(\exp(tH_p)(\xi))=0.
\end{split}\end{equation}
As discussed in \cite{GeMa}, this implies that an interval of
energies around $\lambda$ is non-trapping:
\begin{equation}\begin{split}\label{eq:non-trapping-int}
\exists\delta_0>0\ \text{such that for any}\ \xi\in T^* X^\circ
&\ \text{satisfying}\ p(\xi)\in(\lambda-\delta_0,\lambda+\delta_0),\\
&\ \lim_{t\to\pm\infty}x(\exp(tH_p)(\xi))=0.
\end{split}\end{equation}

The symbolic functional calculus applies in the semiclassical setting as
well -- see \cite{disj} and references given there. 
Here, we will restrict the discussion to the
operator $ P$ given by \eqref{eq:pd}.
Consequently it has
The formula 
\[ f ( P ) = \frac{1}{2 \pi i} \int_{\mathbb C} \bar  \partial_z 
\tilde f ( z ) ( P -z )^{-1} d \bar z \wedge d z \,, \ \ 
\tilde f \in {\mathcal C}_{\mathrm{c}}^\infty ( {\mathbb C}) \,, \
\tilde f|_{\mathbb R} = f \,, \ \bar \partial \tilde f = {\mathcal O}
( |\im z|^\infty ) \,, \]
($ \tilde f $ is an almost analytic extensions of $ f $)
shows
that for $f\in {\mathcal C}^\infty_{\mathrm{c}} (\Real)$, 
$f(P)\in\Psisch^{-\infty ,0,0}(X)$. Also $ 
\sigma^{\star, 0 , 0 }_{h, \mathrm{sc}} ( f ( P ) ) = f ( p ) $.

If $\psi\in\Cinf_c(\Real)$, $\psi\equiv 1$ near $\lambda$, then for $t\in
\Real$,
$1-\psi(\sigma)=
\psit_t(\sigma)(\sigma-(\lambda+it))$, $\psit\in S^{-1}_{cl}(\Real)$
satisfying uniform symbol estimates as $t$ varies over compact sets, so
$\psit(P)\in\Psisch^{-2,0,0}(X)$, and we have proved the following lemma.

\begin{lemma}
Let $P$ be as in \eqref{eq:pd}.
Suppose that $\psi\in\Cinf_c(\Real)$, $\psi\equiv 1$ near $\lambda$, and
suppose that $r,s\in\Real$.
Then there exists $C>0$, independent of $t$ as long as $t$ varies in
compact sets, such that for all $u\in \dist(X)$ with $(P-(\lambda+it))u
\in \Hsc^{r-2,s}(X)$, the following estimate holds:
\begin{equation}\label{eq:ell-3}
\|(\Id-\psi(P))u\|_{\Hsc^{r,s}(X)}
\leq C\|(P-(\lambda+it))u\|_{\Hsc^{r-2,s}(X)}.
\end{equation}
\end{lemma}

\section{Semiclassical estimates}

In this section we will prove the semi-classical resolvent estimates
under the assumption that there exists $q\in S^{0,- \epsilon,0}(X)$, 
$ \epsilon \in (0 , \frac14) $, such that
\begin{gather}
\label{eq:symbol}
\begin{gathered}
2 q H_p q=-b\psi(p)^2 \,, \\  b\in S^{0,1 - 2\epsilon,0}(X) \,, \ \ 
\psi\in\Cinf_c(\Real;[0,1])\,, \ \text{$\psi\equiv 1$ near $\lambda$},
\\ b \geq c_0x^{1 + 2\epsilon }>0 \,.
\end{gathered}
\end{gather}
The existence of $ q $ under global non-trapping assumptions will
be established in Sect.4.

If we write $ Q = {\mathrm{Op}}(q) $ and $ B =
({\mathrm{Op}}(b)+{\mathrm{Op}}(b)^*)/2 $ then,
as reviewed in Sect.2,
\begin{equation}
i[Q^*Q,P]=h \psi(P)B\psi(P)+h^2 R
\end{equation}
with $R\in\Psiscch^{0,2 - 2 \epsilon,0}(X)$. Note that
$\psi(P)\in\Psiscch^{-\infty,0,0}(X)$, i.e.\ it is smoothing, so the
differentiability order $r$ in the weighted Sobolev spaces $\Hsc^{r,s}(X)$
is mostly irrelevant below.
Suppose that $u\in\Hsc^{0, \frac12 + 
\epsilon}(X)$.
Then for $t>0$,
\begin{equation}
\langle u,i[Q^*Q,P]u\rangle=-2\im\langle u,Q^*Q(P-\lambda)u\rangle
-2t\|Qu\|^2.
\end{equation}
(Note that $Q^*QP$ and $PQ^*Q$ are in $\Psiscch^{0,-2\ep,0}(X)$,
so $\langle u,Q^*QP u\rangle$, etc., make sense.)
Thus, taking into account that $2t\|Qu\|^2\geq 0$,
\begin{equation}
h\langle u,\psi(P)B\psi(P)u\rangle
\leq 2|\langle u,Q^*Q(P-(\lambda+it))u\rangle|+h^2 |\langle
u,R u\rangle|.
\end{equation}
By the Cauchy-Schwartz inequality we have, for any $\delta>0$,
\begin{equation}\begin{split}\label{eq:pos-comm-5}
|\langle u,Q^*Q(P-(\lambda+it))u\rangle|&\leq
\|x^{\frac12 + \ep }u\|\|x^{-\frac12 - \ep }Q^*Q(P-(\lambda+it))u\|\\
&\leq
\delta h\|x^{\frac12 + \ep }u\|^2+\delta^{-1}h^{-1}\|x^{-\frac12 - \ep }
Q^*Q(P-(\lambda+it))u\|^2.
\end{split}\end{equation}
Note that $x^{-\frac12 - \ep} R x^{-\frac12 - \ep}\in\Psiscch^{0,1 - 4
\epsilon ,0}(X)$, hence bounded
on $L^2_\scl(X)$ since $ \ep \in (0, \frac14)$. 
Similarly, $x^{-\frac12 - \ep} Q^*Q x^{\frac12 + 3 \epsilon }\in
\Psiscch^{0,0,0}(X)$ is also bounded on the $L^2$ space.
Thus,
\begin{equation}\begin{split}\label{eq:pos-comm-8}
h\langle u,&\psi(P)B\psi(P)u\rangle
-(\delta h+h^2\|x^{{-\frac12 - \ep }}Rx^{-\frac12 - \ep }\|_{\bop(L^2_\scl(X))})
\|x^{{\frac12 + \ep }}u\|^2\\
& \leq \delta^{-1}h^{-1}\|x^{-\frac12 - \ep} Q^*Q x^{
\frac12 + 3 \ep }\|^2_{\bop(L^2_\scl(X))}
\|x^{-\frac12 - 3 \epsilon}(P-(\lambda+it))u\|^2.
\end{split}\end{equation}

We will now use the last assumption in \eqref{eq:symbol}:
$x^{-1 + 2\epsilon}b\psi(p)\geq c_0x^{2\ep }\psi(p)$. Hence by
the sharp G\aa rding estimate,
\begin{equation}
x^{-\frac12 + \ep}\psi(P)B\psi(P)x^{-\frac12 + \ep}
\geq c_0^2x^{2 \ep} \psi(P)^2
x^{2 \ep}+hR_1,\quad R_1
\in\Psisch^{-\infty,1,0}(X).
\end{equation}
Adding $c_0^2 x^{2 \ep}(\Id-\psi(P)^2)x^{2 \ep}$ to both sides gives
\begin{equation}\label{eq:B-norm-5}
x^{{-\frac12 + \ep}}\psi(P)B\psi(P)x^{-\frac12 + \ep}+c_0^2 
x^{2 \ep }(\Id-\psi(P)^2)x^{2 \ep }
\geq c_0^2 x^{4 \ep}+hR_1.
\end{equation}
We also 
note that $|\langle x^{\frac12 - \ep}u,R_1 x^{{\frac12 - \ep}}
u\rangle|\leq C'\|x^{{1- \ep}}u\|^2$.
Thus, applying both sides of \eqref{eq:B-norm-5}
to $x^{{\frac12 - \ep}}u$, and pairing with $x^{\frac12 - \ep} u$
afterwards yields
\begin{equation}\begin{split}
c_0^2\|x^{{\frac12 + \ep}}u\|^2
&\leq
\langle u,\psi(P)B\psi(P)u\rangle
+c_0^2
|\langle(\Id+\psi(P))x^{{\frac12 + \ep}}u,(\Id-\psi(P))x^{{\frac12 + \ep}}u\rangle|
+C' h\|x^{1- \ep }u\|^2\\
&\leq \langle u,\psi(P)B\psi(P)u\rangle+2c_0^2\delta\|x^{\frac12 + \ep}u\|^2
+\delta^{-1} \|(\Id-\psi(P))x^{\frac12 + \ep}u\|^2+C' h\|x^{1 - \ep}u\|^2,
\end{split}\end{equation}
The last term is clearly bounded by $C'h\|x^{\frac12 + \ep}u\|^2$ and 
the second to last term 
can be estimated using \eqref{eq:ell-3}. Choosing
$\delta<1/4$, $h_1=c_0^2/4C'$ gives that for $h\in(0,h_1)$,
\begin{equation}\label{eq:B-norm-2}
\|x^{{\frac12 + \ep}}u\|^2\leq C_1\langle u,\psi(P)B\psi(P)u\rangle
+C_2\|(P-(\lambda+it))u\|^2_{\Hsc^{-2,{-\frac12 - \ep}}(X)}.
\end{equation}
The norm in
the second term
on the right hand side can be replaced by the 
$\Hsc^{-2, \frac12 +  \ep}(X)$ norm.

Combining \eqref{eq:pos-comm-8} and \eqref{eq:B-norm-2},
we thus conclude that there exists $h_0>0$ such that
for $h\in(0,h_0)$,
\begin{equation}
\langle u,\psi(P)B\psi(P)u\rangle
\leq Ch^{-2}\|x^{-\frac12 -3  \ep} (P-(\lambda+it))u\|^2.
\end{equation}
Again using \eqref{eq:B-norm-2}, we conclude that for all $\ep>0$,
\begin{equation}\label{eq:pos-comm-12}
\|u\|_{\Hsc^{0,-1/2-\ep}(X)}
\leq Ch^{-1}\|(P-(\lambda+it))u\|_{\Hsc^{0,1/2+3\ep}(X)},\quad h\in(0,h_0).
\end{equation}
We can modify this argument slightly by inserting $(P+i)(P+i)^{-1}$ in
\eqref{eq:pos-comm-5} between $Q$ and $P-(\lambda+it)$,
to see that the last factor in \eqref{eq:pos-comm-8} can be replaced
by $\|(P+i)^{-1}x^{-\frac12 - 3 \epsilon}(P-(\lambda+it))u\|^2$,
and correspondingly the norm on the right hand side of \eqref{eq:pos-comm-12}
can be replaced by $\|(P-(\lambda+it))u\|_{\Hsc^{-2,1/2+3\ep}(X)}$. A
further slight modification in the same spirit allows us
to conclude that the smoothness
order $r$ in $\Hsc^{r,s}(X)$ can be shifted by the same amount on both sides
of \eqref{eq:pos-comm-12}:
\begin{equation}
\|u\|_{\Hsc^{r,-1/2-\ep}(X)}
\leq Ch^{-1}\|(P-(\lambda+it))u\|_{\Hsc^{r-2,1/2+3\ep}(X)},\quad h\in(0,h_0).
\end{equation}
Now let $u=u_t=R(\lambda+it)f$, $f\in\Hsc^{r,1/2+3\ep}(X)$. Since
$R(\lambda+it)=(P-(\lambda+it))^{-1}\in\Psiscch^{-2,0,0}(X)$
for $t>0$, we see that $u_t\in\Hsc^{r+2,1/2+3\ep}(X)$ for $t>0$.
Thus, the above estimate is applicable and we
conclude that
\begin{equation}\label{eq:unif-res-est}
\|R(\lambda+it)f\|_{\Hsc^{r+2,-1/2-\ep}(X)}
\leq Ch^{-1}\|f\|_{\Hsc^{r,1/2+3\ep}(X)},\quad h\in(0,h_0).
\end{equation}
Note that for a fixed $\psi$, we can let $\lambda$ be arbitrary inside the
region where $\psi\equiv 1$, so a compactness argument gives the uniform
estimate in $\lambda$ as stated in our main Theorem.

\begin{rem}
As in Melrose's paper \cite{mel:spectral},
using these estimates one can show that for fixed $h>0$
the limits
$R(\lambda\pm i0)f$  exist in $\Hsc^{r+2,-1/2-\ep}(X)$ for
$f\in\Hsc^{r,1/2+\ep}(X)$, $\ep>0$. Hence, \eqref{eq:unif-res-est} yields
\begin{equation}
\|R(\lambda+i0)f\|_{\Hsc^{r+2,-1/2-\ep}(X)}
\leq Ch^{-1}\|f\|_{\Hsc^{r,1/2+\ep}(X)},\quad h\in(0,h_0),
\end{equation}
as well.
\end{rem}

\section{Symbol construction}
Let $p$ be the principal symbol of $P$. Thus, near $\bX$,
\begin{equation}
p=\tau^2+g_\pa(y,\mu)+x^\gamma r,\quad r\in S^{2,0}(X),
\end{equation}
where $g_\pa$ is the metric on the boundary, and we denote the metric
function on the cotangent bundle the same way.
Its Hamilton vector field $H_p$
is of the form
\begin{equation}
x(2\tau (x\pa_x+\mu\cdot\pa_\mu)-2g_\pa \pa_\tau+H_{g_\pa})+x^{ 1 + \gamma }
W,
\quad W\in\Vb(\sct X) \otimes S^{1,0} ( X) \,.
\end{equation}
Here we will be mainly concerned with the $(x,\tau)$ variables, so we
rewrite this as
\begin{equation}
H_p=x(2\tau+x^\gamma a) (x\pa_x)-x(2g_\pa+x^\gamma b) \pa_\tau +2x\tau \mu\cdot\pa_\mu
+xH_{g_\pa}+x^{1+\gamma}W',
\end{equation}
where $a,b\in S^{1,0}( X)$, and
$W'$ is now a vector field tangent to the $\bX$ fibers, i.e.\ it is
a vector field in $\pa_y$ and $\pa_\mu$ with coefficients in 
$ S^{1,0} ( X) $. In this section we take
$\lambda^2$, not $\lambda$, as the spectral parameter.

\begin{figure}[ht]
\begin{center}
\mbox{\epsfig{file=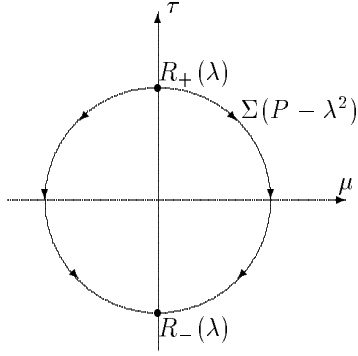}}
\end{center}
\caption{The projection of the characteristic variety $\Sigma(P-\lambda^2)$
and the bicharacteristics of $H_p$ inside it to the $(\tau,\mu)$-plane.}
\label{fig:bich}
\end{figure}

As indicated,
we make the assumption that a small interval of energies around
$\lambda^2$ is non-trapping, i.e.
\begin{equation}\begin{split}\label{eq:non-trapping-2}
\exists\delta_0>0\ \text{such that for any}\ \xi\in T^* X^\circ
&\ \text{satisfying}\ p(\xi)\in(\lambda^2-\delta_0,\lambda^2+\delta_0),\\
&\ \lim_{t\to\pm\infty}x(\exp(tH_p)(\xi))=0.
\end{split}\end{equation}

Now,
\begin{equation}
H_p (x^{-1}\tau)=-2(\tau^2+g_\pa)+x^\gamma f,\ f\in S^{1,0} (X)\,,
\end{equation}
so there exists $\ep_1>0$ such that
for $\xi\in\sct X$ satisfying $p(\xi)\in(\lambda^2/2,2\lambda^2)$,
$x<\ep_1$, $-(H_p (x^{-1}\tau))(\xi)\geq c_0>0$. Since $p$ is
constant along integral curves of $H_p$, we see that if
$\exp(-tH_p)(\xi)$, $t\geq T$,
stays in $x<\ep_1$ (which holds under our non-trapping assumption for
sufficiently large $T$),
then $x^{-1}\tau$ tends to $+\infty$; in particular
$\tau$ is non-negative for all large $t$.
By reducing $\ep_1>0$ if necessary, we also see that
there exist $\delta_1>0$, $\ep_1>0$ such that
for $\xi\in\sct X$,
\begin{equation}
|p(\xi)-\lambda^2|<\delta_1,\ x(\xi)<\ep_1,\ |\tau|<7\lambda/8
\Rightarrow g_\pa(\xi)\geq c_1>0.
\end{equation}
Reducing
$\ep_1>0$ further if necessary, we can thus arrange that
\begin{equation}
|p(\xi)-\lambda^2|<\delta_1,\ x<\ep_1,\ |\tau|<7\lambda/8
\Rightarrow -x^{-1}H_p\tau(\xi)\geq c_1>0.
\end{equation}
Thus, we see that given any $x_0>0$, $\xi\in T^*X^\circ$ with
$|p(\xi)-\lambda^2|<\delta_1$,
there exists $T>0$ such that
\begin{equation}\label{eq:non-tr-tau}
t\geq T\Rightarrow
\tau(\exp(-tH_p)(\xi))>2\lambda/3,\ x(\exp(-tH_p)(\xi))<x_0/2.
\end{equation}

We now define a symbol $q\in S^{-\infty,0}(\sct X)$ whose most important
properties are that
\begin{equation}
q\geq 0\Mand x^{-1}H_p q\leq 0.
\end{equation}
We will always use a localization
in the energy via a factor $\psi(p)$ where $\psi\in\Cinf_c(\Real)$ is
supported in $(\lambda^2-\delta,\lambda^2+\delta)$, where $\delta\in(0,
\lambda^2)$ is
a fixed small constant with $\delta<\delta_1$, $\delta_1$ as above.
Let
\begin{equation}
M=\sup \{|a(\xi)|+|b(\xi)|:\ p(\xi)\leq 2\lambda^2\}<+\infty;
\end{equation}
here we used that
$p^{-1}((-\infty,2\lambda^2])$ is a compact subset of $\sct X$.
Also, let
\begin{equation}
x_0=\min\{(\lambda/6(M+1))^{1/\gamma},(c_1/2(M+1))^{1/\gamma},
\ep_1\}.
\end{equation}
Let $\chi_-\in\Cinf(\Real)$ be supported in $(\lambda/3,+\infty)$,
identically $1$ on $(2\lambda/3,+\infty)$, with $\chi'_-\geq 0$,
and similarly let
$\chi_+\in\Cinf(\Real)$ be supported in $(-\infty,-\lambda/3)$,
identically $1$ on $(-\infty,-2\lambda/3)$, with $\chi'_+\leq 0$.
Also, let $\chi_\pa\in\Cinf_c(\Real)$ be supported in
$(-7\lambda/8,7\lambda/8)$, with $\chi'_\pa\geq (6\lambda/c_1)\chi_\pa\geq 0$
on $(-7\lambda/8,3\lambda/4)$, and $\chi_\pa(-3\lambda/4)>0$.
Let $\rho\in\Cinf_c([0,+\infty))$ be identically $1$ on $[0,1/2]$, supported
in $[0,1)$, $\rho'\leq 0$ on $[0,+\infty)$.
In the incoming region
we will take the symbol
\begin{equation}
q_-=x^{-\epsilon}\chi_-(\tau)\psi(p)\rho(x/x_0),
\end{equation}
in the outgoing one the symbol
\begin{equation}
q_+=x^{\ep}\chi_+(\tau)\psi(p)\rho(x/x_0),
\end{equation}
with $ \ep \in (0, \frac14) $.
In the intermediate region we take
\begin{equation}
q_\pa=x^{-\ep }\chi_\pa(\tau)\psi(p)\rho(x/x_0).
\end{equation}

Note that for any $\alpha\in\Real$, $\chi,\rho\in\Cinf(\Real)$,
\begin{equation}\begin{split}
&x^{-\alpha-1}H_p(x^\alpha\chi(\tau)\rho(x/r))\\
&=(2\tau+x^\gamma a)(\alpha\rho(x/r)+
r^{-1}\rho'(x/r))\chi(\tau)
-(2g_\pa(\xi)+x^\gamma b)\rho(x/r)\chi'(\tau).
\end{split}\end{equation}
Note that in the definition of $q_-$, $\alpha=- \ep<0$, so $\alpha\rho(x/r)+
r^{-1}\rho'(x/r)\leq 0$ everywhere.
Moreover, on $\supp\chi_-$,
$\tau>\lambda/3>0$, so for $x(\xi)\leq x_0$,
$\xi\in\supp\psi(p)$, $\tau(\xi)\in\supp\chi_-$,
$2\tau+x^\gamma a\geq \lambda/3>0$. In addition, $\tau\leq 2\lambda/3$ on
$\supp\chi'_-$, so if $\xi\in\supp(\rho(x/x_0)\chi'(\tau)\psi(p))$ then
$g_\pa\geq c_1>0$, hence $2g_\pa+x^\gamma b\geq c_1>0$ there.
Thus,
\begin{equation}\label{eq:H_p-q_--gl}
x^{-1 + \ep}H_p q_-\leq 0.
\end{equation}
Moreover, $x\leq x_0/2$ implies $\rho'(x/x_0)=0$, and $\tau\geq 2\lambda/3$
implies $\chi'_-(\tau)=0$, so
\begin{equation}\label{eq:H_p-q_--loc}
x\leq x_0/2,\ \tau\geq 2\lambda/3\Rightarrow -x^{-1 + \ep }
H_p q_-\geq c_2\psi(p),
\ c_2>0.
\end{equation}

The difference between $q_-$ and
$q_+$ is that $\tau\rho'$ is positive on $\supp
\chi_+$, and $-\chi'_+$ is also positive,
so the negativity estimate only holds away from $\supp\rho'$ and
$\supp\chi'_+$.
Thus, there is no analogue of \eqref{eq:H_p-q_--gl}, but the following
analogue of \eqref{eq:H_p-q_--loc} still holds:
\begin{equation}
x\leq x_0/2,\ \tau\leq -2\lambda/3\Rightarrow -x^{-1-\ep}H_p q_+\geq c_3\psi(p),
\ c_3>0.
\end{equation}

Next, $q_\pa$ provides the connection between the incoming and outgoing
regions.
Since $\chi_\pa'$ can be used to estimate $\chi_\pa$ on $(-7\lambda/8,
3\lambda/4)$, we see that
\begin{equation}
\tau(\xi)\in(-7\lambda/8,3\lambda/4),\ x(\xi)\leq x_0/2,\ \xi\in
\supp\psi(p)\Rightarrow
|(2\tau+x^\gamma a)\chi_\pa(\tau)|\leq c_1\chi'_\pa(\tau)/2.
\end{equation}
Since $\alpha=- \ep$, $|\alpha|<1$, so we conclude that
\begin{equation}
\tau(\xi)\in(-7\lambda/8,3\lambda/4),\ x(\xi)\leq x_0/2,\ \xi\in
\supp\psi(p)\Rightarrow
-x^{-1 + \ep }H_p q_\pa\geq c_1\chi'_\pa(\tau)\psi(p)\geq 0.
\end{equation}
Note that on $(-3\lambda/4,3\lambda/4)$, $\chi'_\pa\geq C>0$, so
\begin{equation}
\tau(\xi)\in(-3\lambda/4,3\lambda/4),\ x(\xi)\leq x_0/2,\ \xi\in
\supp\psi(p)\Rightarrow
-x^{-1 + \ep }H_p q_\pa\geq c_4\psi(p),\ c_4>0.
\end{equation}

For $\xi\in T^*X^\circ$
with $p(\xi)\in(\lambda^2-\delta_0,\lambda^2+\delta_0)$,
take $T=T_\xi>0$ as in \eqref{eq:non-tr-tau}, so for
$t\geq T$ we have
$\tau(\exp(-tH_p)(\xi))>2\lambda/3$, $x(\exp(-tH_p)(\xi))<x_0/2$.
We will define a symbol $q_\xi$ which is supported in a neighborhood
of the bicharacteristic segment $\{\exp(-t H_p)(\xi):\ t\in[0,T+1]\}$, and
which satisfies $H_p q\leq 0$ over
\begin{equation}
K'=\{\xi'\in T^*X^\circ:\ x(\xi')\geq x_0/2\Mor (x(\xi')\leq x_0/2\Mand
\tau(\xi')\leq 2\lambda/3\}.
\end{equation}
Namely, let $\Sigma$ be a hypersurface through $\xi$ which is transversal
to $H_p$. Then there is a neighborhood $U_\xi$ of $\xi$, such that
$V_\xi=\{\exp(-t (U_\xi\cap\Sigma)):\ t\in(-1,T+2)\}$ is a neighborhood of
the above bicharacteristic segment, which we can think of as a product
$(-1,T+2)\times(U_\xi\cap\Sigma)$, and $(T+1/2,T+2)\times(U_\xi\cap\Sigma)$
is disjoint from $K'$. Now let $\phi_\xi\in\Cinf_c(U_\xi\cap\Sigma)$
be identically $1$ near $\xi$, and let $\chi_\xi\in\Cinf_c(\Real)$
be supported in $(-1,T+2)$, $\chi_\xi\geq 0$, $\chi'_\xi\geq 0$ on
$(-1,T+2/3)$.
Using the product coordinates, we can think of $\phi_\xi$ and $\chi_\xi$
as functions of $\sct X$ with compact support in $V_\xi$. Let
\begin{equation}
q_\xi=\chi_\xi \phi_\xi \psi(p),
\end{equation}
so
\begin{equation}
H_p q_\xi=-\chi'_\xi \phi_\xi \psi(p).
\end{equation}
Thus, for $\xi'\in K'$, $H_p q_\xi(\xi')\leq 0$.

Now let $K\subset T^*X^\circ$ be the compact set
\begin{equation}
K=\{\xi\in T^*X^\circ:\ \xi\in\supp\psi(p),\ x(\xi)\geq x_0/4\}.
\end{equation}
Since $K$ is compact, applying the previous argument for
every $\xi\in K$ gives a $U_\xi$, and a $U'_\xi\subset U_\xi$ on which
$\phi_\xi=1$. Since $\{U'_\xi:\ \xi\in K\}$ covers $K$, the compactness
of $K$ shows that we can pass to a finite subcover, $\{U'_{\xi_j}:
\ j=1,\ldots,N\}$. We let
\begin{equation}
q_\circ=\sum_{j=1}^N q_{\xi_j}.
\end{equation}

\begin{figure}[ht]
\begin{center}
\mbox{\epsfig{file=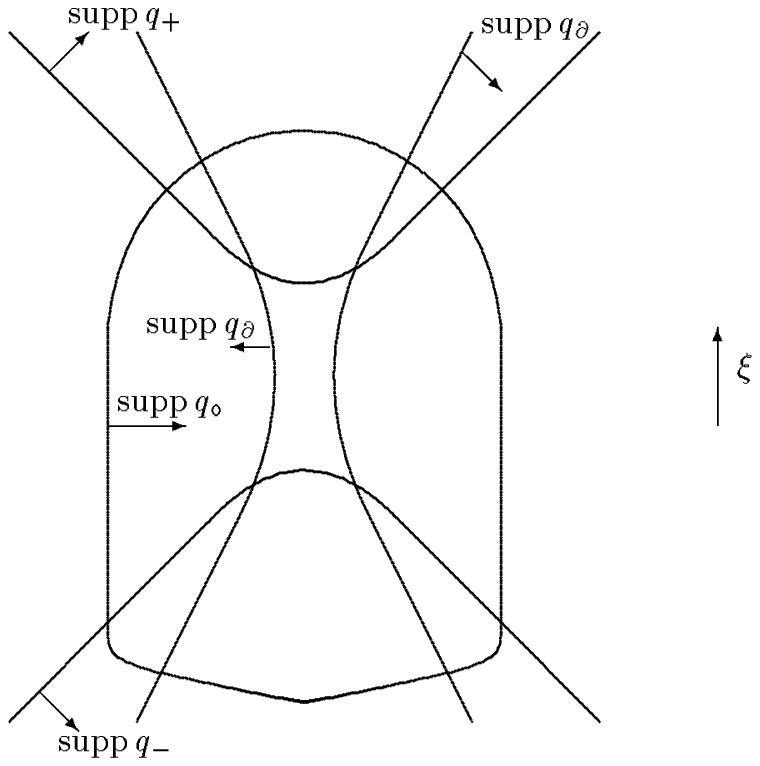}}
\end{center}
\caption{Supports of $q_+$, $q_-$, $q_\circ$ and $q_\partial$ for
$X=\overline{\Rn}$. $\sct X$ is identified with
$\overline{\Rn}\times\Rn_\xi$, and the covector $\xi$ is fixed on the picture.}
\label{fig:supports}
\end{figure}

The symbol we use in the positive commutator estimate is
\begin{equation}
q=q_-+C''q_\pa+C q_\circ+C'q_+,
\end{equation}
with $C,C',C''>0$ chosen appropriately. Namely note that in the region
$x\leq x_0/2$, $\tau\geq 2\lambda/3$, which is the only place where
$H_pq_\circ$ is positive,
we have the estimate $-x^{-1 + \ep }H_p q_-\geq c_2>0$.
Since $x^{-1 + \ep }H_pq_\circ$ is bounded, we can choose $C>0$ sufficiently
small so that $-x^{-1 + \ep }H_p (q_-+Cq_\circ)$
is still bounded below by a positive
constant in this region. Then $-x^{-1 + \ep }
H_p (q_-+Cq_\circ)$ is non-negative
everywhere, and it is bounded
below by a positive constant on $x\geq x_0/2$ as well as on $x\leq x_0/2$,
$\tau\geq 2\lambda/3$.
But this is the only region where the bounded function
$x^{-1 + \ep }H_p q_\pa$ is positive, so by choosing
$C''>0$ sufficiently small, we can arrange that
$-x^{-1 + \ep }
H_p (q_-+Cq_\circ+C''q_\pa)$ is non-negative everywhere, and it is
bounded below by a positive constant on $x\geq x_0/2$, as well as on
$x\leq x_0/2$, $\tau\geq -3\lambda/4$.
But this is the only region where $x^{-1-\ep}H_p q_+>0$. Thus,
by choosing $C'>0$ sufficiently small, and taking into account that
$x^{-1 + \ep }H_p=x^{2 \ep }x^{-1 - \ep }H_p q_+$, with $x^{-1-\ep}H_p q_+$ as well as
$x^{2 \ep}$ bounded, we can arrange that
$-x^{-1+\ep}H_p q$ is non-negative everywhere, and
$-x^{-1-\ep}H_p q$ bounded below by a positive constant everywhere.
In summary, we have proved the proposition needed in Sect.3 
(see \eqref{eq:symbol})

\begin{prop}
There exist functions $q\in S^{-\ep,\infty}(X)$, $\psi\in\Cinf_c(\Real;
[0,1])$, $\psi\equiv 1$ near $\lambda^2$, and $c',c''>0$ such that
\begin{equation}
q\geq c'x^{\ep}\psi(p),\ -H_p q\geq c''x^{1+\ep}\psi(p).
\end{equation}
\end{prop}

Thus, the results of the previous section show that
there exists $h_0>0$ such that
for $h\in(0,h_0)$,
\begin{equation}
\|R(\lambda^2+it)f\|_{\Hsc^{*,-1/2-\ep}(X)}\leq C_0h^{-1}\|f\|
_{\Hsc^{*,1/2+3\ep}(X)}.
\end{equation}


\begin{thebibliography}{XX}

\bibitem{Ag} Agmon, Sh., Spectral theory of Schr\"odinger operators on 
Euclidean and non-Euclidean spaces, {\em Comm. Pure Appl. Math.} 
{\bf 39} (1986), S3-S16.

\bibitem{BrPe} Bruneau, V., and Petkov, V., Semiclassical resolvent 
estimates for trapping perturbations, preprint, 1999.

\bibitem{disj}
  Dimassi, M., and Sj\"ostrand,~J., Spectral asymptotics in the
  semi-classical limit, Lecture Notes to appear.

\bibitem{Ge} G\'erard, Ch., Semiclassical resolvent estimates for two
and three body Schr\"odinger operators, {\em Comm. P.D.E.}, {\bf 15}
(1990), 1161-1178.

\bibitem{GeMa} G\'erard, Ch., and Martinez, A., Principe d'absorption 
limite pour des op\'erateurs de Schr\"odinger 
\`a longue port\'ees, {\em C.R. Acad. Sci. Paris} {\bf 306} (1988), 121-123.

\bibitem{HelSj} Helffer, B., and  Sj\"ostrand,~J. 
Resonances en limite semi-classique. {\em M\'em. Soc. Math. France (N.S.) }
\textbf{24-25} (1986).

\bibitem{Ho:e} H\"ormander, L., On the existence and the regularity of 
solutions of linear pseudo-differential equations. 
{\em Enseignement Math. } {\bf 17}(2) (1971), 99-163. 

\bibitem{ho:v3} H\"ormander, L., \emph{Linear partial differential
    equations}, vol.3, Springer Verlag, Berlin, 1985.


\bibitem{JeMoPe} Jensen, A., Mourre, E. and Perry, P., Multiple commutator
estimates and resolvent smoothness in quantum scattering theory, 
{\em Ann. Inst. H. Poincar\'e (phys.~th\'eor.)} {\bf 41} (1984), 207-225.


\bibitem{mel:spectral} Melrose, R.~B., {Spectral and scattering theory for
    the Laplacian on asymptotically Euclidean spaces,} \emph{Spectral and
    scattering theory} (M.~Ikawa, ed.), Marcel Dekker, 1994, 85--130.

\bibitem{MelZw} Melrose, R.~B., and Zworski, M., 
Scattering metrics and geodesic flow at infinity. {\em Invent. Math.}
\textbf{124}(1996), 389-436.

\bibitem{RoTa} Robert, D., and Tamura, H., Semiclassical estimates for
resolvents and asymptotics for total scattering cross-sections, 
{\em Ann. Inst. H. Poincar\'e (phys.~th\'eor.)} {\bf 47} (1987), 415-442.
 
\bibitem{Wa} Wang, X.P., Time-decay of scattering solutions and classical
trajectories, 
{\em Ann. Inst. H. Poincar\'e (phys.~th\'eor.)} {\bf 47} (1987), 25-37.

 
\bibitem{WZ1} Wunsch, J., and Zworski, M., 
Distribution of resonances for asymptotically Euclidean manifolds,
preprint, 1999.

\end{thebibliography}
\end{document}